\documentclass[10pt]{amsart}
\usepackage{amsmath}
\usepackage{amssymb,amsfonts,amscd}
\usepackage{graphicx}

\begin{document}

\title{Simplices with equiareal faces}
\author{Victor Alexandrov}
\address{Sobolev Institute of Mathematics, Koptyug ave., 4, Novosibirsk, 630090, Russia and
Department of Physics, Novosibirsk State University, Pirogov str., 2, Novosibirsk, 630090, Russia}
\email{alex@math.nsc.ru}
\author{Nadezhda Alexandrova}
\address{Institute of Mining, Krasnyj pr., 54, Novosibirsk, 630091, Russia and
Department of Physics, Novosibirsk State University, Pirogov str., 2, Novosibirsk, 630090, Russia}
\email{alex@math.nsc.ru}
\author{Gunter Wei\ss}
\address{Technische Universit\"at Dresden, Fachrichtung Mathematik,
Institut f\"ur Geometrie, 01062 Dresden, Germany}
\email{Gunter.Weiss@tu-dresden.de}
\subjclass{51M20}
\keywords{Heron formula, area, tetrahedron, pyramid, simplex}
\date{\today}
\begin{abstract}
We study simplices with equiareal faces in the Euclidean 3-space
by means of elementary geometry. We present an unexpectedly simple
proof of the fact that, if such a simplex is non-degenerate, than 
every two of its faces are congruent. We show also that 
this statement is wrong for degenerate simplices and find
all  degenerate simplices with equiareal faces.
\end{abstract}
\maketitle

\renewcommand{\thefootnote}{\fnsymbol{footnote}}

\section{Introduction}\label{section1}
\footnotetext{The first author was supported in part by the Russian State Program for Leading Scientific Schools, Grant~NSh--8526.2008.1.}

This paper deals with  simplices with equiareal faces in the Euclidean 3-space.
A simplex is called a \textit{simplex with equiareal faces} if all its faces 
have the same area and is called \textit{degenerate} if all its vertices lie
on a single plane. Our primary interest is the following 

\textbf{Problem 1.} 
\textit{Prove that all faces of any non-degenerate 
simplex with equiareal faces in the Euclidean 3-space
are congruent to each other.}

Problem 1 already appeared in several contexts.
We mention here some of those appearances though we cannot 
say that our knowledge is complete.

In the 1960s Professor Hans Vogler in Vienna used Problem 1 
to convince his students how powerful the synthetic geometry is.
As far as we know, his
desciptive-geometrical solution was never published.

Independently, in the 1970's and 1980's 
the following form of Problem 1 was used  
to prevent `undesirable applicants'
from joining the Moscow State University:
`The faces of a triangular pyramid have the same area. 
Show that they are congruent.'
An elementary solution to the latter problem
given in \cite{Vardi} shows that indeed the problem is 
far from being trivial.

In \cite{McMullen}, among other results P. McMullen has
proved that \textit{a non-degenerate simplex is 
a simplex with equiareal faces
if an only if its opposite 
pairs of edges have  the same lengths.}
It is easy to see that this statement is equivalent to
Problem 1. 
McMullen's proof is very short and natural, but it is not
elementary since it rests on Minkowski's theorem 
on uniqueness and existence of a~closed convex polyhedron
with given directions
and areas of faces (see, e.\,g., \cite{Alexandrov};
the direction of a~face is determined by
the outward unit normal to the face).

In 2007 Professor Robert Connelly in Cornell University
has brought our attention to 
the fact that it is reasonable to study degenerate simplices
with equiareal faces. He argued that when we fix three vertices,
say $A$, $B$, and $C$, and move arbitrarily the fourth vertex,
say $D$, nothing special happens when $D$ occurs 
in the plane~$ABC$: the degenerate simplex is obtained as a limit
of non-degenerate ones; the notions of the vertex, edge, and
face are clear for it; the notion of 
the face area is well defined; from combinatorial point of
view there is no difference between degenerate and 
non-degenerate simplices. Besides, the degenerate
polyhedra play very important role in `advanced' study
of convex polytopes, see, e.\,g., \cite{Alexandrov}.

Also Professor Robert Connelly has brought our attention to 
the fact that, in contrast with Problem 1, there are 
degenerate simplices with non-congruent equiareal faces.
In fact, every parallelogram equipped with its diagonals may be 
treated as a degenerate simplex with equiareal faces as well as
every four points on a line may be treated as a vertex-set 
of a degenerate simplex with equiareal faces 
(with all areas equal zero). The problem is whether
there are some other degenerate simplices with equiareal faces.

In Section~\ref{section2} we give the shortest available for us
elementary solution to Problem~1. 
As far as we know, its main idea should be 
attributed to Professor Hans Vogler
who is now at the University of Innsbruck.
In Section~\ref{section3} we use  Heron's formula to study all
simplices with equiareal faces, both degenerate 
and non-degenerate.

\section{A short elementary solution to Problem 1}\label{section2}

Note that in order to solve Problem 1 it is sufficient to solve the
following problem~2 which is of independent interest.

\textbf{Problem~2.} 
\textit{Let $ABCD$ be a non-degenerate simplex, let 
$a(\triangle ABC){=}a(\triangle ABD)$, and let
$a(\triangle ACD)=a(\triangle BCD)$, where 
$a(\triangle XYZ)$ stands for the area of the 
triangle $\triangle XYZ$.
Prove that $|AC|=|BD|$ and $|BC|=|AD|$, where
$|XY|$ stands for the length of the straight 
line segment~$XY$} (see Fig.~1).

\textbf{Solution} to Problem 2.
Let $P$ be the plane which is parallel to the line~$AB$ 
and contains the line~$CD$ (see Fig.~2). 
We are going to study the orthogonal projection of 
the simplex $ABCD$ into the plane $P$. 
Denote by~$X^\perp$ the image of a point $X$ under that projection
and denote by~$X^\ast$ the foot of the perpendicular to the line~$AB$
emanated from the point~$X$.  
\begin{figure}
\includegraphics[width=0.45\textwidth]{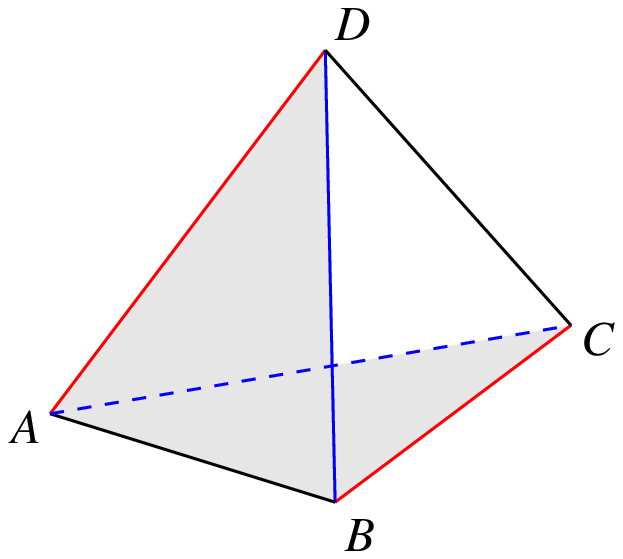}\hfill
\includegraphics[width=0.54\textwidth]{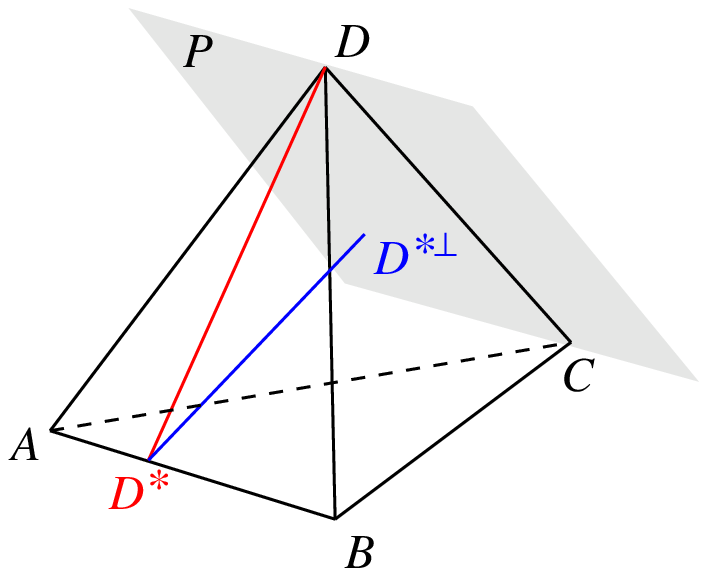}\\
\hskip-65mm\parbox[t]{0.47\textwidth}{\caption{}}\hskip10mm
\parbox[t]{0.47\textwidth}{\caption{}}
\end{figure}

Since~$AB$ is parallel to~$P$ we 
have~$|D^\ast D^{\ast\perp}|=|C^\ast C^{\ast\perp}|$.  
Since the triangles~$\triangle ABC$ and~$\triangle ABD$ 
have equal areas and the common side~$AB$ we conclude 
that~$|DD^\ast |=|CC^\ast |$. 
Applying Pythagoras theorem to the right-angled
triangles~$\triangle DD^\ast D^{\ast\perp}$ and 
$\triangle CC^\ast C^{\ast\perp}$ 
we get~$|DD^{\ast\perp}|^2
=|DD^{\ast\perp}|^2-|D^\ast D^{\ast\perp}|^2
=|CC^{\ast\perp}|^2-|C^\ast C^{\ast\perp}|^2
=|CC^{\ast\perp}|^2$.  
In terms of the quadrilateral  $A^\perp CB^\perp D$ this means that 
the vertices~$D$ and~$C$ lie at the same distance from 
the diagonal~$A^\perp B^\perp$.  
Note also that the triangles $\triangle A^\perp B^\perp D$ and 
$\triangle A^\perp B^\perp C$ have the same area and the points~$D$ 
and~$C$ lie on the different sides of the line passing
through the points~$A^\perp$ and~$B^\perp$. 
In fact, if they lie on the same side, the line through
the points~$A^\perp$ and~$B^\perp$ must be parallel to the line 
through the points~$D$ and~$C$ and, thus,
the  points~$A$, $B$, $C$, and~$D$ should be coplanar. 
A contradiction.

Similar arguments applied to the triangles~$\triangle ACD$
and~$\triangle BCD$ show that the vertices~$A^\perp$ and~$B^\perp$ 
of the quadrilateral~$A^\perp CB^\perp D$ lie at the same distance 
from the line through the points~$D$ and~$C$ and, moreover, 
lie on the different sides of that line. Hence, 
the quadrilateral~$A^\perp CB^\perp D$ is convex. 

Recall that if a convex planar quadrilateral is such that 
every two opposite vertices lie at the same distance from  
the diagonal that joins the rest two vertices then the 
quadrilateral is a parallelogram.  This implies that
the quadrilateral~$A^\perp CB^\perp D$ is a parallelogram 
and, thus, that $|A^\perp C|=|B^\perp D|$  and 
$|B^\perp C|=|A^\perp D|$.

Applying Pythagoras theorem to the right-angled triangles 
$\triangle AA^\perp C$ and $\triangle BB^\perp D$ and 
taking into account that $|AA^\perp|=|BB^\perp|$
we get $|AC|^2=|AA^\perp|^2+|A^\perp C|^2
=|BB^\perp|^2+|B^\perp D|^2=|BD|^2$.

Similar arguments applied to the triangles
$\triangle AA^\perp D$   and $\triangle BB^\perp C$  
show that  $|BC|=|AD|$. \hfill Q.E.D.
                                    
\section{A study of degenerate and non-degenerate simplices\\ with equiareal faces}\label{section3}

From Section \ref{section1} 
we know the following three types of equiareal 
non-degenerate and degenerate simplices in the Euclidean 3-space:

Type 1: non-degenerate simplex with all faces congruent;

Type 2: parallelogram equipped with its diagonals; and

Type 3: four points on a line treated as a vertex-set 
of a degenerate simplex with equiareal faces of zero area.

In this Section we use Heron's formula to study the following

\textbf{Problem~3.} 
\textit{Prove that there are no simplices with equiareal faces which do not belong to Types 1--3.}

\textbf{Solution} to Problem 3.
Let $T$ be a simplex with equiareal faces. Let the first face 
(I) of $T$ has edges of the lengths $a$, $b$, and $c$; 
the second face 
(II) has edges of the lengths $a$, $y$, and $z$; 
the third face 
(III) --- $b$, $x$, and $z$; and the forth face 
(IV) --- $c$, $x$, and $y$. 
(Equivalently, we can say that side $x$ is opposite to $a$; 
side $y$ --- to $b$; and $z$ --- to $c$). 
Let $S$ be common area of the faces of $T$.

Heron's formula for the face (I) yields
\begin{align}
(4S)^2 &=(a+b+c)(-a+b+c)(a-b+c)(a+b-c)\notag\\
 &= 2a^2b^2+2a^2c^2+2b^2c^2-a^4-b^4-c^4=-(a^2-b^2+c^2)^2+4a^2c^2.
\end{align}

Now let's use Heron's formula to express the fact that the faces 
(I) and (II)  have the same areas
\begin{align}
(a+y+z)(-& a+y+z)(a-y+z)(a+y-z)-(4S)^2\notag\\
=& 2a^2y^2+2a^2z^2+2z^2y^2-a^4-y^4-z^4-(4S)^2\\
=&-(z^2-y^2-a^2)^2-(4S)^2+4y^2a^2=0.\notag
\end{align}
Solving this equation with respect to $z^2$ yields
\begin{equation}
z^2=y^2+a^2\pm\sqrt{4a^2y^2-(4S)^2}.
\end{equation}

Similarly, we use Heron's formula to express the fact that 
the faces (I) and (III)  have the same areas
\begin{align}
(b+z+x)(-& b+z+x)(b-z+x)(b+z-x)-(4S)^2\notag\\
=& 2b^2z^2+2b^2x^2+2z^2x^2-b^4-z^4-x^4-(4S)^2\notag\\
=& -(z^2-x^2-b^2)^2-(4S)^2+4b^2x^2=0.\notag
\end{align}
Solving this equation with respect to $z^2$ yields
\begin{equation}
z^2=x^2+b^2\pm\sqrt{4b^2x^2-(4S)^2}.
\end{equation}

At last, we use Heron's formula to express the fact that the faces 
(I) and (IV)  have the same areas
\begin{align}
(c+x+y)(-& c+x+y)(c-x+y)(c+x-y)-(4S)^2\notag\\
=& 2c^2x^2+2c^2y^2+2x^2y^2-c^4-x^4-y^4-(4S)^2\notag\\
=& -(y^2-x^2-c^2)^2-(4S)^2+4c^2x^2=0.\notag
\end{align}
Solving this equation with respect to $y^2$ 
\begin{equation}
y^2=x^2+c^2\pm\sqrt{4c^2x^2-(4S)^2}.
\end{equation}

Eliminate $z^2$ from (3) and (4)
$$
y^2+a^2\pm\sqrt{4a^2y^2-(4S)^2}=x^2+b^2\pm\sqrt{4b^2x^2-(4S)^2},
$$
then twice square this equation in order to eliminate square 
roots and use the
formula $(a^2-b^2+c^2)^2=4a^2c^2-(4S)^2$ (which is a consequence 
of (1)) to obtain
\begin{gather}
 4a^4y^4+4b^4x^4+(x^2+b^2-y^2-a^2)^4-8a^2b^2x^2y^2-2a^2y^2(x^2+b^2-y^2-a^2)^2\notag\\
 -2b^2x^2(x^2-y^2-a^2+b^2)^2=-4(x^2-y^2-a^2+b^2)^2.
\end{gather}

So, we arrive at the most computationally difficult, but 
still straightforward
point of the solution: substitute $y^2$ in (6) by the 
right-hand side of (5).
After simplifications and multiple usage of the formula  
$(a^2-b^2+c^2)^2=4a^2c^2-(4S)^2$ we get
\begin{equation}
4(c^2x^2-4S^2)(x^2-a^2)^2 S^4= \bigl[x^4-x^2(a^2+b^2-c^2)-a^2(b^2-c^2)^2\bigr]S^4.
\end{equation}

We see that $S=0$ is a root of (7) which corresponds to simplices of Type 3. 
In order to find the other roots,
cancel $S^4$ in (7), rearrange terms and use the formula  
$(a^2-b^2+c^2)^2=4a^2c^2-(4S)^2$ again to arrive at
$$
(x^2-a^2)^2\bigl[x^4-2x^2(b^2+c^2)+a^2(2b^2+2c^2-a^2)\bigr]=0.
$$
Note that the bi-quadratic expression in the brackets has two roots: 
$x^2=a^2$ and $x^2=2b^2+2c^2-a^2$.
Note also that we can obtain similar equations for $y$ and 
$z$ just by permuting three pairs of symbols $(a,x)$, $(b,y)$, 
and $(z,c)$.
As a result we find 8 solutions for $x, y,$ and $z$ that 
are accumulated as rows in the following table:

\vskip10pt
\begin{center}
\renewcommand{\arraystretch}{1.5}
\begin{tabular}{|c|c|c|c|}
\hline
 & $x$ & $y$ & $z$ \\
\hline
Solution 1 & $a$ & $b$ & $c$ \\
\hline
Solution 2 & $a$ & $b$ & ${\sqrt{2a^2+2b^2-c^2}}^{}$ \\
\hline
Solution 3 & $a$ & $\sqrt{2a^2-b^2+2c^2}$ & $c$ \\
\hline
Solution 4 & $a$ &  $\sqrt{2a^2-b^2+2c^2}$  & $\sqrt{2a^2+2b^2-c^2}$ \\
\hline
Solution 5 & $\sqrt{2b^2+2c^2-a^2}$ & $b$ & $c$ \\
\hline
Solution 6 & $\sqrt{2b^2+2c^2-a^2}$ & $b$ & $\sqrt{2a^2+2b^2-c^2}$  \\
\hline
Solution 7 & $\sqrt{2b^2+2c^2-a^2}$ & $\sqrt{2a^2-b^2+2c^2}$  & $c$ \\
\hline
Solution 8 & $\sqrt{2b^2+2c^2-a^2}$ & $\sqrt{2a^2-b^2+2c^2}$  & $\sqrt{2a^2+2b^2-c^2}$ \\
\hline
\end{tabular}
\end{center}
\vskip10pt

Solution 1, obviously, corresponds to the simplices of Type 1. 

Solutions 2, 3, and 5 correspond to simplices of Type 2.
For example, a simplex $T$, corresponding to Solution 2,  
is degenerated into a parallelogram with the side lengths $a$ 
and $b$ and the diagonals $c$ and $\sqrt{2a^2+2b^2-c^2}$ 
(recall that in any parallelogram the sum of the squared 
lengths of all sides equals the sum of the squared lengths 
of the both diagonals).

Solutions 4, 6, and 7 correspond to simplices of Type 2 again, but 
the face (I), with edge lengths $a$, $b$, and $c$, must be a 
right-angled triangle this time. For example, consider 
Solution 4. We have 
\begin{equation}
x^2=a^2,\quad y^2=2a^2-b^2+2c^2, \quad\mbox{and}\quad z^2=2a^2+2b^2-c^2. 
\end{equation}
Using the formula (2) we get $-(z^2-y^2-a^2)^2-(4S)^2+4a^2y^2=0$.  
Now we use the formula (1) and, after some simplifications, 
we get $a^4-(b^2-c^2)^2=0$. Without loss of generality, 
we may assume that $b\geqslant c$. This yields $a^2+c^2=b^2$ 
and, thus, the face (I) is a right-angled triangle. 
Moreover, now the formula (8) implies that $y^2= b^2$ and 
$z^2= 4a^2+c^2$. Hence, the simplex $T$, corresponding to 
Solution 4,  is degenerated to a parallelogram with the side 
lengths $a$, $b$, $x=a$, and $y=b$ and the diagonals of the 
lengths $c$ and $z=\sqrt{4a^2+c^2}=\sqrt{2a^2+2b^2-c^2}$. 
Hence, the simplex $T$ is of Type 2. Solutions 6 and 7 are 
treated similarly.

Solution 8 does not correspond to any simplex in the Euclidean 
3-space (neither degenerated nor non-degenerated). 
In fact, we have 
\begin{equation}
x^2=2b^2+2c^2-a^2, \quad y^2=2a^2-b^2+2c^2, \quad\mbox{and}\quad z^2=2a^2+2b^2-c^2.
\end{equation}
Using the formula (2) we get $-(z^2-y^2-a^2)^2-(4S)^2+4a^2y^2=0$.  
Now we use formula (1) and, after some simplifications, we get 
$a^4-(b^2-c^2)^2=0$ or 
\begin{equation}
(a^2-b^2+c^2)(a^2+b^2-c^2)=0.
\end{equation}
The geometric meaning of the formula (10) is that either $b$ or $c$ is the hypotenuse 
of the right-angled triangle (I) with the sides $a$, $b$, and $c$.
Similarly we can substitute (9) into formulas (3) and (4). 
Proceeding as above we get
\begin{gather}
(a^2+b^2-c^2)(-a^2+b^2+c^2)=0,\\
(a^2-b^2+c^2)(-a^2+b^2+c^2)=0.
\end{gather}
The geometric meaning of the formula (11) is that either $c$ or $a$  
is the hypotenuse of the right-angled triangle (I) with the sides 
$a$, $b$, and $c$.
Similarly, the formula (12) implies that either $b$ or $a$  is 
the hypotenuse of the right-angled triangle (I) with the sides 
$a$, $b$, and $c$.
But the triangle (I) has only one hypotenuse! Hence the equations 
(10)--(12) can not hold true simultaneously. 
This means that Solution 8
does not correspond to any simplex.

Now we can conclude that there is no 
simplices with equiareal faces which do not belong to Types 1--3.
\hfill{Q.E.D.}

\end{document}